%% file: comparison.tex
\documentclass[12pt]{amsart}
\usepackage[dvips]{graphicx}
\usepackage{color}

\def\RR{\mathbb R}
\def\HH{\mathbb H}
\def\ts{{\mathcal T}(S)}
\def\a{\alpha}
\def\s{\sigma}
\def\L{{\Lambda}}

\def\T{{\mathcal T}}

\def\cs{{{\mathcal C}(S)}}
\def\mul{\stackrel{.}{\asymp}}
\def\add{ \stackrel{\mbox{\tiny $+$}}{\asymp} }
\def\precadd{ \stackrel{\mbox{\tiny $+$}}{\prec} }
\def\precmul{ \stackrel{.}{\prec} }
\def\succmul{ \stackrel{.}{\succ} }
\def\dd{\partial}
\def\r{\rho}
\def\ep{\epsilon}
\def\Teich{Teichm\"uller\ }

\newcommand{\st}
{%
  \hspace{0.3em}\left|\rule{0mm}{2.0ex}\right.\hspace{0.3em}
}
\newcommand{\set}[1]{\left\{\hspace{0.3em}#1\hspace{0.3em}\right\}}

\newtheorem{theorem}{Theorem}[section]

\newtheorem{lemma}[theorem]{Lemma}
\newtheorem{proposition}[theorem]{Proposition}
\newtheorem{corollary}[theorem]{Corollary}
\newtheorem{remark}[theorem]{Remark}

\newtheorem{introthm}{Theorem}
\newtheorem{introcor}[introthm]{Corollary}

\title{Comparison Between Teichm\"uller and Lipschitz Metrics}
\author{Young-Eun Choi} 
\email{choiye@math.ucdavis.edu}

\author{Kasra Rafi}
\email{rafi@math.uconn.edu}
\urladdr{http://www.math.uconn.edu/~rafi/} 

\begin{document}
\maketitle

\section{Introduction}
\noindent
The Teichm\"uller distance between two points $\s,\tau$ in \Teich
space $\T(S)$ is defined in terms of the minimal quasiconformal
constant $K(\s,\tau)$ between $\s$ and $\tau$.  In \cite{thurston}
Thurston introduced an analogous metric on $\T(S)$ by considering the
least possible value of the global Lipschitz constant
$\Lambda(\s,\tau)$ from $\s$ to $\tau$.  On the one hand, Kerckhoff
\cite{kerckhoff} showed that $K(\s,\tau)$ can be formulated in terms
of the ratio of extremal lengths of simple closed curves:
\begin{equation}
\label{eqn:kerckhoff}
K(\s,\tau)=\sup_\alpha
\frac{\mbox{Ext}_{\tau}(\alpha)}{\mbox{Ext}_{\s}(\alpha)}
\end{equation}
and on the other, it was shown by Thurston \cite{thurston} that the
minimal Lipschitz constant $\Lambda(\s,\tau)$ is given by the ratio of
lengths in the hyperbolic metric:
\begin{equation}
\label{eqn:thurstonratio}
\Lambda(\s,\tau)=\sup_\alpha
\frac{l_{\tau}(\alpha)}{l_{\sigma}(\alpha)}.
\end{equation}

A comparison of $K(\s,\tau)$ and the ratio of lengths in
 Equation~\eqref{eqn:thurstonratio} was first given by 
Wolpert \cite{wolpert}, who proved that for any $K$--quasiconformal
map $f$ from $\s$ to $\tau$ and any simple closed curve $\alpha$
$$\frac{l_\tau(f(\alpha))}{l_\s(\alpha)} \leq K.$$ This implies,
in particular, that
\begin{equation*}
\Lambda(\s,\tau) \leq K(\s,\tau).
\end{equation*}

In this paper, we compare the \Teich and Lipschitz metrics by
comparing the two ratios in Equations~\eqref{eqn:kerckhoff} and
\eqref{eqn:thurstonratio}.  Our method is to analyze the ratio of
hyperbolic lengths in much the same way the ratio of extremal lengths
was analyzed by Minsky in \cite{minskyprod} for the purpose of showing
that certain regions in the thin part of Teichm\"uller space have
product structures.  However, since $K(\s,\tau)$ is symmetric and
$\L(\s,\tau)$ is not \cite{thurston}, it is necessary to choose some
symmetric version of $\L$ to make the comparison more meaningful.
Thus, we take
$$L(\s,\tau)=\max\{\L(\s,\tau),\L(\tau,\s)\}$$
and define the \Teich and Lipschitz metrics,
respectively, as follows
\begin{align*}
d_\T(\s,\tau) &= \frac{1}{2} \log K(\s,\tau) \\
d_L(\s,\tau) &= \log L(\s,\tau). 
\end{align*}
Note that the factor of $1/2$ has been left out in the Lipschitz
metric. This is due to the fact that on the thick part of Teichm\"uller
space, we can compare the two metrics up to an additive error, as we
shall shortly see.

Although $\L(\s,\tau)$ is not symmetric, it is easy to check that it
satisfies the following ordered triangle inequality:
$$\L(\rho,\tau) \leq \L(\rho,\s) + \L(\s,\tau)$$
and further,
satisfies the property that $\L(\s,\tau)=0$ if and only if $\s=\tau$.
Thus $d_L(\s,\tau)$ defines a genuine metric, in
that it is symmetric, takes the value zero if and only if $\s = \tau$,
and satisfies the triangle inequality.  In \cite{papa}, it was shown
that on the \Teich space of the torus, the \Teich metric and a
similarly defined Lipschitz metric are, in fact, equal. In
contrast, we show that for a hyperbolic surface $S$, the two metrics
are not comparable. In particular,

\begin{introthm}
\label{thm:example}
There are sequences $\s_n,\tau_n
\in \T(S)$ such that, as $n \rightarrow \infty$,
$$d_L(\s_n,\tau_n) \rightarrow 0,
 \ \ d_\T(\s_n,\tau_n) \rightarrow \infty.$$
\end{introthm}

As is often the case, however, no incongruities occur on the thick
part of \Teich space, and there the two metrics are quasi-isometric to
one another.  In fact, they are equal up to a bounded additive error.
This is a consequence of the following theorem proved in
Section~\ref{sec:thick}:
\begin{introthm}
\label{thm:thick0}
 For $\rho \in \T(S)$, let $\mu_\rho$ be
a short marking for $\rho$.
For every $\ep>0$, there is a constant $c$ depending on $\ep$ such
that, for any $\s,\tau$ in the $\ep$--thick part of $\T(S)$, the
following quantities differ from one another by at most $c$:
\begin{align*}
1.& \ d_\T(\sigma,\tau)&    2.& \ d_L(\sigma,\tau) \\
3.& \ \log \, \max_{\alpha \in \mu_\s} 
    \frac{l_{\tau}(\alpha)}{l_{\sigma}(\alpha)} &
4.& \ \log \, \max_{\alpha \in \mu_\tau}
     \frac{l_{\sigma}(\alpha)}{l_{\tau}(\alpha)}. 
\end{align*}
\end{introthm}
In particular, in order to estimate the \Teich distance between two
points in the thick part, one need only compare the lengths of a
finite number of curves with respect to the two metrics. 

To compare the metrics on the thin part of \Teich space, we prove in
Section~\ref{sec:prodreg} an analogue of Minsky's product region
theorem \cite{minskyprod}.  Let $\Gamma$ be a collection of $k$
disjoint, homotopically distinct, simple closed curves on $S$ and let
$Thin_\ep(S,\Gamma)$ be the set of $\s \in \T(S)$ such that
$l_\s(\gamma) \leq \ep$ for all $\gamma \in \Gamma$.  Let
$T_\Gamma=\T(S \setminus \Gamma) \times U_1 \times \cdots \times U_k$,
where $S \setminus \Gamma$ is the analytically finite surface obtained
from $S$ by pinching all the curves in $\Gamma$ and where $U_i$ is the
subset $\{ (x,y): y\geq 1/\ep \}$ of the upper-half plane.  The
Fenchel-Nielsen coordinates on $\T(S)$ give rise to a natural
homeomorphism $\Pi:Thin_\ep(S,\Gamma) \rightarrow \T_\Gamma$.
 Then Minsky's product region theorem states:
\begin{theorem}[Minsky \cite{minskyprod}]
\label{thm:minskyprod}
 Let $d_{\T_\Gamma}$ be the sup metric
$$d_{\T_\Gamma}=\sup \Big\{d_{\T(S \setminus \Gamma) },
\frac{1}{2}d_{\HH_1}, \ldots, \frac{1}{2} d_{\HH_k} \Big\}$$
on $\T_\Gamma$, where $d_{\T(S \setminus \Gamma)}$ is the \Teich metric
on $\T(S\setminus \Gamma)$ and $d_{\HH_i}$ is the restriction of the
hyperbolic metric on the upper-half plane to $U_i$.  Then, for $\ep$
sufficiently small, there is a constant $c$ depending on $\ep$, such
that for any $\s,\tau \in Thin_\ep(S,\Gamma)$,
$$|d_\T(\s,\tau) - d_{\T_\Gamma}(\Pi(\s),\Pi(\tau))| < c.$$
\end{theorem}
In the analogue for
the Lipschitz metric, we define the sup metric
$$d_{L_\Gamma}=\sup \{ d_{L(S \setminus \Gamma)}, d_{L(A_1)},\ldots,
d_{L(A_k)} \}$$
on $\T_\Gamma$, where $d_{L(S \setminus \Gamma)}$ is
the Lipschitz metric on $\T(S\setminus \Gamma)$ and $d_{L(A_i)}$ is a
modification of the hyperbolic metric on $U_i$ (see
Section~\ref{sec:prodreg} for details):

\begin{introthm} 
\label{introthm:split}
  For $\ep$ sufficiently small, there is a constant $c$ depending on
  $\ep$, such that for any $\s,\tau \in Thin_\ep(S,\Gamma)$,
$$|d_L(\s,\tau) - d_{L_\Gamma}(\Pi(\s),\Pi(\tau))| < c.$$
\end{introthm}
A more precise statement is given in Theorem~\ref{thm:prodreg}.  Our
proof follows parallel to Minsky's, but requires only elementary
hyperbolic geometry, since we need not deal with extremal lengths.

As a consequence of Theorem~\ref{thm:thick0}, one can deduce the following 
purely combinatorial result. For a subsurface $Z$, let $d_Z(\mu_1,\mu_2)$ 
be the distance between the projections of $\mu_1$ and $\mu_2$ to $Z$, 
measured in the arc complex of $Z$ (see \cite{MM}, \cite{rafi2} for details).

\begin{introcor}
\label{cor:intro} 
There is a constant $k$ such that  for any markings $\mu_1$ and 
$\mu_2$ on $S$,
\begin{equation} \label{eqn:intersection}
\log i(\mu_1,\mu_2) \asymp \sum_{Y}
\left[ d_Y(\mu_1,\mu_2) \right]_k+
\sum_{A} \log \left[ d_A (\mu_1,\mu_2) \right]_k,
\end{equation}
where $Y$ ranges over all subsurfaces of $S$ that are not annuli,
$A$ ranges over all annuli, and where $[x]_k = 0$ if $x<k$ and 
$[x]_k =x$ if $x \geq k$.
\end{introcor}

In \cite{MM}, Masur and Minsky provide an estimate, similar to the
right-hand side of \eqref{eqn:intersection}, for the number of
elementary moves needed to change $\mu_1$ to $\mu_2$. Using their
result and examining how the intersection number between two markings
changes as a result of applying a sequence of elementary moves to one
of them, one can show that the right-hand side of \eqref{eqn:intersection} 
is an upper bound for $\log i(\mu_1,\mu_2)$
(there is no clear combinatorial argument for proving the inequality
in the other direction). In this context, Corollary~\ref{cor:intro}
states that, along an efficient path in the marking space, the
intersection number increases at the fastest possible rate.

\subsection{Notation}
Often, we shall compare two functions $f,g$ on $\T(S)$ and use the
notation $f \prec g$, $f \asymp g$ to mean, respectively, that there
are positive constants $k,c$ such that $f \leq k g + c$, $\frac{1}{k}
g - c \leq f \leq k g + c$. We also use $f \precmul g$,
$f \add  g$ to mean, respectively, that there is only a
multiplicative constant, or only an additive constant, involved.  In
particular, $f \mul 1$ means that the function $f$ is bounded both
above and below by positive constants.  The constants $k$ and $c$ usually
depend on the topological type of $S$, which will not be subsequently
mentioned. Any other dependencies will be explicitly noted.

\section{The thick part} \label{sec:thick}
\noindent
Let $S$ be a surface of finite topological type.  Given $\ep>0$, the
$\ep$--thick part of \Teich space is the set of $\s \in \T(S)$ such
that the infimum of the injectivity radius measured in $\s$, taken
over all points in $S$, is greater than $\ep$. When we simply say
``the thick part'', we mean it is the $\ep$--thick part for some $\ep$
which has already been chosen.

A {\em marking} on $S$ is a collection of homotopically distinct,
simple closed curves in $S$ obtained by first choosing a pants curves
system, i.e., a collection of mutually disjoint curves that cut $S$
into pairs of pants (where a hole may be a puncture of $S$) and then
by choosing an additional collection of curves that together with the
pants system cuts the surface into disks and punctured disks. To make
the choice of a marking less arbitrary, additional conditions on the
choice of curves are often specified.

For $ \s \in \T(S)$, we define a {\em short marking} $\mu_\s$, as
follows.  First choose a pants system by taking the shortest curve in
$S$, then the next shortest curve disjoint from the first, and so on
until a complete pants system $\underline{\a}$ is formed. We remark
that throughout this paper, when we say the ``length of a curve'', we
always mean the length of its geodesic representative.  Next, choose a
``dual'' curve $\delta_\a$ for each $\a \in \underline{\a}$ that is
disjoint from $\underline{\a} \setminus \a$, and that is shortest
among all such curves. There may be a finite number of possible short
markings for $\s$.

A lemma of Bers' says that there is a uniform constant $N$ such that
every $\s \in \T(S)$ has a pants curves system $\underline{\a}$ with
the property that $l_\s(\a) < N$ for all $\a \in \underline{\a}$.
Hence, if $\s$ is in the $\ep$--thick part of $\T(S)$ so that all the
curves in a short marking $\mu$ have length bounded below as well,
then the lengths of the dual curves are bounded above and so
$l_\s(\mu)$ is bounded above by some quantity depending only on $\ep$.
Conversely, given a marking $\mu$ and a number $B>0$, the metrics $\s
\in \T(S)$ such that $l_\s(\mu)=\sum_{\a \in \mu}l_\s(\a) \leq B$ has
bounded diameter in $\T(S)$, where the bound depends only on $B$ (see
for example \cite{minskytop}).  Thus there is a coarse correspondence
between the thick part of \Teich space and the set of markings. This
idea is implicit in the theorems to follow.

\setcounter{introthm}{1}
\begin{introthm} 
\label{thm:thick} 
For every $\ep>0$, there is a constant $c$ depending on $\ep$ such
that, for any $\s,\tau$ in the $\ep$--thick part of $\T(S)$, the
following quantities differ from one another by at most $c$:
\begin{align*}
1.& \ d_\T(\sigma,\tau)&    2.& \ d_L(\sigma,\tau) \\
3.& \ \log \, \max_{\alpha \in \mu_\s} 
    \frac{l_{\tau}(\alpha)}{l_{\sigma}(\alpha)} &
4.& \ \log \, \max_{\alpha \in \mu_\tau}
     \frac{l_{\sigma}(\alpha)}{l_{\tau}(\alpha)}. 
\end{align*}
\end{introthm}
First we need the following lemma. Let $g:\RR \to \T(S)$ be the \Teich
geodesic that passes through $\s$ and $\tau$ and let $q_t$ be the
family of quadratic differentials representing $g$.  We assume all
quadratic differential metrics have been normalized to have area $1$.

\begin{lemma}
\label{lem:exponential}
For every marking $\mu$ on $S$ there exist $l_0$ and $t_0$ such that 
$$l_{q_t}(\mu) \mul l_0 \, e^{|t-t_0|}.$$
\end{lemma}
\begin{proof}
Recall that a quadratic differential $q_t$ defines a pair of
measured foliations on the surface $S$, called the horizontal and
the vertical foliations. For every curve $\a$ the horizontal length
$h_t(\a)$ of $\a$ is the intersection number of $\alpha$ with the
horizontal foliation and the vertical length $v_t(\a)$ of $\a$ is
the intersection number of $\a$ with the vertical foliation. We have
(see for example \cite{rafi1})
$$l_{q_t}(\a) \mul h_t(\a) + v_t(\a). $$
Let $t_\a$ be the time when $\a$ is balanced, i.e., the time when 
the horizontal length and the vertical length of $\a$ are equal.
Let $l_\a = l_{q_{t_\a}}(\a)$. Along a Teichm\"uller 
geodesic, the horizontal length of $\a$ increases and the vertical 
length of $\a$ decreases exponentially fast. Therefore,
$$l_{q_t}(\alpha) \mul l_\alpha \cosh(t-t_\a).$$ 
Thus, for every marking $\mu$
\begin{equation} \label{eqn:f}
l_{q_t}(\mu)=\sum_{\a \in \mu} l_{q_t}(\a) \mul
\sum_{\a \in \mu} l_\a\cosh (t-t_\a).
\end{equation}
Denote the right hand side of \eqref{eqn:f} by $f(t)$. Let $t_0$ be the 
time when $f(t)$ is minimum and let $l_0 = f(t_0)$. Since 
$$\cosh(t-t_\a) \leq \cosh (t_0-t_\a) \, e^{|t-t_0|} ,$$ 
we have
\begin{equation}
\label{eqn:upper}
\sum_{\a \in \mu} l_\a\cosh (t-t_\a) \leq
\sum_{\a \in \mu} l_\a \cosh (t_0 - t_\a) \, e^{|t-t_0|} = 
l_0 \, e^{|t-t_0|}.
\end{equation}

To prove the inequality in the other direction, we observe that
the derivative of $f(t)$ with respect to $t$ at $t=t_0$
is $\sum_\a l_\a \sinh(t_0-t_\a)=0$,
which implies
$$
\sum_{\a \in \mu} l_\a e^{t_0-t_\a}=\sum_{\a \in \mu} l_\a
e^{t_\a-t_0}=l_0.
$$
If $n$ is the number of curves in $\mu$, the above equation implies
that  there exist
$\beta,\gamma \in \mu$ such that
$$
l_{\beta} \, e^{t_0 - t_\beta} \geq \frac{l_0}{n}
\quad \text{and} \quad l_{\gamma} \, e^{t_\gamma - t_0} \geq \frac{l_0}{n}.
$$
Thus we have
\begin{align}
f(t) 
 &=\sum_{\a \in \mu} l_\a\cosh (t-t_\a) \geq
  l_\beta \cosh (t-t_\beta) + 
  l_\gamma \cosh (t-t_\gamma)\nonumber\\
 &\geq
  \frac{1}{2}\Big[l_\beta \, e^{t-t_0} e^{t_0 - t_\beta} +
  l_\gamma \, e^{t_\gamma-t_0} e^{t_0 - t}\Big] \nonumber\\
\label{eqn:lower}
&\geq \frac{l_0}{2n} \, e^{|t-t_0|}.
\end{align}
Equations~\eqref{eqn:upper} and \eqref{eqn:lower} show that
$f(t) \mul l_0 \, e^{|t-t_0|}$. This and \eqref{eqn:f} prove 
the lemma.
\end{proof}

\begin{proof}[Proof of Theorem~\ref{thm:thick}]
We show that the first three quantities are comparable, the proof
for the remaining term is similar.  Suppose that for $a<b$, we have
$g(a)=\s$, $g(b)=\tau$ so that $d_\T(\s,\tau)=b-a$.  Since the
moduli space of the thick part is compact, we know that the
hyperbolic lengths of curves in $\s,\tau$ are proportional to their
quadratic differential lengths in $q_a, q_b$, respectively (see
\cite{rafi3} for a more general discussion).  Therefore, there are
multiplicative constants depending only on $\ep$, such that for any
simple closed curve $\a$,
\begin{equation}
\label{eqn:qlengths}
\frac{l_\tau(\a)}{l_\s(\a)}
\mul \frac{l_{q_b}(\a)}{l_{q_a}(\a)}.
\end{equation}
Moreover, since
$$
  l_{q_b}(\a) \mul l_\a \cosh (b-t_\a) 
  \leq e^{b-a} l_\a \cosh (a-t_\a)= e^{b-a} l_{q_a}(\a),
$$
it follows from Equation~\eqref{eqn:qlengths} that
$$d_L(\s,\tau) \precadd d_\T(\s,\tau).$$
Thus it remains to be shown that there is a curve $\a \in \mu_\s$
such that
$$ b-a \precadd \log \frac{l_{q_b}(\alpha)}{l_{q_a}(\alpha)}. $$

Let $l_{q_t}(\mu_\s) \mul l_0 \, e^{|t-t_0|}$ as in
Lemma~\ref{lem:exponential}.  Then,
\begin{equation}
\label{eqn:q_b}
 l_0 \, e^{|b-a| - |a-t_0|} \precmul 
 l_{q_b}(\mu_\s) \precmul  l_0 \, e^{|b-a| + |a-t_0|}.
\end{equation}
First we show that $|a-t_0|$ is bounded above.  Since $\s$ is in the
thick part, the $q_a$--length and the $\sigma$--length of $\mu_\sigma$
are comparable to one another. Moreover, since $\mu_\s$ is a short
marking in $\s$, its $\s$--length is bounded both above and below.
Therefore, we have:
\begin{equation}
\label{eqn:l_0}
l_0 \, e^{|a-t_0|} \mul l_{q_a}(\mu_\s) \mul l_{\s}(\mu_\s) \mul 1.
\end{equation}
Furthermore, we can see that $l_0$ is bounded below, as follows.  A
marking divides the surface into disks and punctured disks.  For any
quadratic differential $q$, the $q$--area of a disk or a punctured disk
is less than the square of its perimeter.  Therefore,
we have for all $t$,
\begin{equation}
\label{eqn:q-area} 
1=\mbox{area}_{q_t}(S) \stackrel{.}{\prec}
 \sum_{\alpha \in \mu_\s} l_{q_t}(\alpha)^2.
\end{equation}
Applied to $t=t_0$ we get $l_0 \stackrel{.}{\succ} 1$. It then follows
from Equation~\eqref{eqn:l_0} that $|a-t_0| \stackrel{.}{\prec} 1$, as
desired.  Thus, it follows from Equation~\eqref{eqn:q_b} that
$$l_{q_b}(\mu_\s) \mul e^{b-a}.$$
But, as we saw in Equation~(\ref{eqn:l_0}), since the
$q_a$--lengths of curves in $\mu_\s$ are bounded above and below, it
follows that there exists a curve $\a \in \mu_\s$ such that
$$l_{q_b}(\a) \mul l_{q_a}(\a) e^{b-a},$$
which is what we wanted.
\end{proof}

\begin{theorem} 
\label{thm:d=i}
Let $\sigma$ and $\tau$ be points in the $\ep$--thick part of
Teichm\"uller space and let $\mu_\s$ and $\mu_\tau$ be their short
markings, respectively.  Then there is an additive constant depending
only on $\ep$, such that
$$d_\T(\sigma, \tau) \add \log \, i(\mu_\s,\mu_\tau),$$
where
$i(\mu_\s,\mu_\tau)$ is the total number of intersections between the
curves in $\mu_\s$ and the curves in $\mu_\tau$.
\end{theorem}
\begin{proof} The $\tau$--length of a curve is proportional
to its intersection number with $\mu_\tau$ (see for example
\cite[Lemma 4.7]{minskytop}). Therefore,
\begin{equation}
\label{eqn:intersection1}
i(\mu_\s,\mu_\tau) \mul \sum_{\alpha \in \mu_\s} l_{\tau}(\alpha)
\mul  \max_{\alpha \in \mu_\s} l_{\tau}(\alpha).
\end{equation}
Since $\s$ is in the thick part of $\T(S)$, we have 
$l_\s(\a) \mul 1$ for every curve $\a \in \mu_\s$. Thus, it follows
from Theorem~\ref{thm:thick} that
\begin{equation}
\label{eqn:intersection2}
\log \max_{\alpha \in \mu_\s} l_{\tau}(\alpha) \add
\log \max_{\alpha \in \mu_\s}
\frac{l_{\tau}(\alpha)}{l_{\sigma}(\alpha)} \add 
d_\T(\sigma, \tau).
\end{equation}
The theorem follows from Equations~\eqref{eqn:intersection1} 
and \eqref{eqn:intersection2}.
\end{proof}

\begin{remark}
The above theorem implies that the logarithm of the intersection number is 
almost a distance function on the marking space.  In particular, it
satisfies a quasi-triangle inequality. That is, for markings $\mu_1$, $\mu_2$
and $\mu_3$ we have
$$\log i(\mu_1, \mu_3) \precadd \log i(\mu_1, \mu_2) + \log i(\mu_2,
\mu_3).$$
This ``distance function'' is similar, but not comparable to
the distance defined on the space of markings in \cite{MM}.
\end{remark}

\begin{proof}[Proof of Corollary~\ref{cor:intro}] 
For given markings $\mu_1$ and $\mu_2$, one can find points $\s_1$
and $\s_2$ in the thick part of \Teich space such that $\mu_1$ and
$\mu_2$ are short markings in $\s_1$ and $\s_2$, respectively. In
\cite{rafi2}, a combinatorial formula is given for the \Teich
distance between two points in the thick part of \Teich space.  It
states that $d_\T(\s_1, \s_2)$ is comparable to the right-hand side
of Equation~\eqref{eqn:intersection}. Also, Theorem~\ref{thm:d=i}
states that $\log i(\mu_1, \mu_2) \add d_\T(\s_1, \s_2)$.  Together
these two results prove the corollary.
\end{proof}

\section{Product regions in the Lipschitz metric} \label{sec:prodreg}
\noindent
In this section, we prove the analogue of Minsky's product region
theorem for the Lipschitz metric. 
\subsection{An $(\ep_0,\ep_1)$--decomposition}
First, we need to recall the notion of an
$(\ep_0,\ep_1)$--decomposition defined in \cite{minskyprod}.  Let $0 <
\ep_1 < \ep_0$ be two numbers less than the Margulis constant.  Let
$\s$ be a hyperbolic metric on $S$ and suppose
$\gamma_1,\ldots,\gamma_k$ are geodesics with length $l_\s(\gamma_i)
\leq \ep_1$. Let $A_1,\ldots,A_k$ be the collection of annular
neighborhoods of $\gamma_1,\ldots,\gamma_k$, respectively, such that
the boundary components of $A_i$ each have length $\ep_0$.  
A component $Q$ of $S \setminus \cup A_i$ is called a {\em hyperbolic} 
component and the
entire collection $\mathcal P$ of hyperbolic components and annular
components is called an $(\ep_0,\ep_1)$--{\em decomposition}.  We
assume that $\ep_0,\ep_1$ are chosen so that any simple geodesic that
intersects an annular component $A$ is either the core of $A$ or is
made up of arcs that run from one boundary component of $A$ to
another.  We remark that in \cite{minskyprod}, what we have described
is called a {\em partial} $(\ep_0,\ep_1)$--decomposition. There, the
term $(\ep_0,\ep_1)$--decomposition is reserved for the case where
$\{\gamma_1,\ldots,\gamma_k\}$ is the full set of curves whose length
$l_\s$ satisfies $l_\s \leq \ep_1$.

In the course of arguments to follow, we shall further require that
$\ep_0 /\ep_1 >2$ so that certain desired estimates hold (see for
example Lemma~\ref{lem:arctocurve}).  We therefore assume
$\ep_0,\ep_1$ have been chosen once and for all to satisfy all the
conditions stated above and henceforth use the notation $f \mul g$, $f
\add g$, etc., to mean that the multiplicative or additive constants
which appear depend only on this choice of $\ep_0,\ep_1$ (and on the
topological type of $S$).

\subsection{Decomposing the length of a curve}
\label{sec:decomposition}
Consider the intersection of a simple closed curve $\zeta$ with the
components of an $(\ep_0,\ep_1)$--decomposition.  For a hyperbolic
component $Q$, let ${\mathcal C}(Q,\dd Q)$ denote the homotopy classes
of simple closed curves in $Q$ and of essential arcs in $Q$ with
endpoints on $\dd Q$, under homotopies that keep any endpoints of arcs
on $\dd Q$.  Define the orthogonal projection $\zeta_Q$ of $\zeta$ to
be the geodesic representative of $\zeta \cap \,Q$ in ${\mathcal
  C}(Q,\dd Q)$ that has the shortest length (see \cite[\S
2.3]{minskyprod}).  In particular, every arc in $\zeta_Q$ is
perpendicular to $\dd Q$. It is not hard to show the following:
\begin{proposition}
\label{prop:thickthin}
Let ${\mathcal P}$ be the components of an $(\ep_0,\ep_1)$--decomposition
for $\s$ and let $Q,A \in {\mathcal P}$ be respectively, a hyperbolic and
annular component.  Then, for any simple closed curve $\zeta$, the
following estimates hold:
\begin{align}
\label{eqn:Qlength0}
i(\zeta,\dd Q) & \succmul \left| l_\s(\zeta \cap Q)- l_\s(\zeta_Q)
\right| \\ \label{eqn:Alength0} i(\zeta, \gamma) & \succmul \left|
  l_\s(\zeta \cap A) - \Big[\log \frac{\epsilon_0}{l_\s(\gamma)} +
  l_\s(\gamma) \, \frac {Tw_\s(\zeta,\gamma)}2 \Big] i(\zeta, \gamma)
\right|
\end{align}
where $\gamma$ is the core geodesic of $A$.
\end{proposition}
Here, $Tw_\s(\zeta,\gamma)$ is the absolute value of the twist of
$\zeta$ around $\gamma$ defined in \cite[\S 3]{minskyprod}. In
Equation~\eqref{eqn:Alength0}, $\log [\ep_0/l_\s(\gamma)]$ is
approximately half the width of $A$; the right-hand side describes the
sum of lengths of piecewise geodesic arcs homotopic to $\zeta \cap A$
in ${\mathcal C}(A,\dd A)$, each of which goes perpendicularly from one
component of $A$ to $\gamma$, wraps around $\gamma$ a number of
$Tw_\s(\zeta,\gamma)$ times (up to an error of $1$), then goes out the
other end of $A$ orthogonally. The idea is that most of the twisting
$\zeta$ does around $\gamma$, takes place in $A$ \cite{minskyprod}.
This is also the reason that Equation~\eqref{eqn:Qlength0} is true.
For a proof, see \cite{CRS}.

Since the components of $\dd Q$ each have a collar of some definite
width, $l_\s(\zeta \cap Q) \succmul i(\zeta,\dd Q)$ and $l_\s(\zeta_Q)
\succmul i(\zeta, \dd Q) $.  Similarly, since $\gamma$ has a collar of
definite width, terms in the righthand side of
Equation~\eqref{eqn:Alength0} are larger than a multiple of $i(\zeta,
\gamma)$. Therefore, we can rewrite Equations~\eqref{eqn:Qlength0} and
\eqref{eqn:Alength0} as follows:
\begin{corollary}
\label{cor:thickthin}
Let $Q,A$ be as in Proposition~\ref{prop:thickthin}. Then for
any simple closed curve $\zeta$ on $S$, we have
\begin{align}
\label{eqn:Qlength}
l_\s(\zeta \cap Q) &\mul l_\s(\zeta_Q) \\
\label{eqn:Alength}
l_\s(\zeta \cap A)  &\mul \left [\log \frac{\epsilon_0}{l_\s(\gamma)}
 + l_\s(\gamma)\cdot \frac{Tw_\s(\zeta,\gamma)}2 \right] i(\zeta, \gamma).
\end{align}
\end{corollary}

\subsection{Regular Annuli} 
Let $A$ be an annulus. We call a metric $\rho$ on $A$ a 
{\em regular metric} if $(A, \rho)$ is isometric to quotient of some closed
neighborhood $ \set {p \in \HH^2 \st d_{\HH^2}(p, G) \leq r }$ of a
geodesic $G$ in $\HH^2$, by a hyperbolic isometry with axis $G$.  For
$\ep>0$, let $U_\ep(A)$ be the space of all regular metrics on $A$
such that the core of $A$ has length at most $\ep$ and such that each
component of $\dd A$ has length $\ep_0$. Two metrics are considered
equivalent if they differ by an isotopy of $A$ fixing $\dd A$
pointwise. Define the distance  between $\rho_1,\rho_2 \in U_\ep(A)$ to be
$$
  d_{L(A)}(\rho_1,\rho_2)=\sup_{\beta \in {\mathcal C}(A,\dd A)}  
  \Bigg| \log \frac{l_{\rho_1}(\beta)}{l_{\rho_2}(\beta)} \Bigg|,
$$
where ${\mathcal C}(A,\dd A)$ is the set of homotopy classes of
non-trivial simple loops or arcs in $A$, under homotopies that fix
the endpoints.  As usual, the length $l_\rho(\beta)$ means the
length of the $\rho$--geodesic representative of $\beta$.
Clearly $d_{L(A)}(\rho_1,\rho_2)$ is symmetric, and is zero if and only
if $\rho_1=\rho_2$.  To see that the triangle inequality holds,
observe that
\begin{align*}
\Bigg| \log \frac{l_{\rho_1}(\beta)}{l_{\rho_2}(\beta)}\Bigg|+
\Bigg| \log \frac{l_{\rho_2}(\beta)}{l_{\rho_3}(\beta)}\Bigg| &
\geq \log \frac{l_{\rho_1}(\beta)}{l_{\rho_2}(\beta)} +
\log \frac{l_{\rho_2}(\beta)}{l_{\rho_3}(\beta)}=
\log \frac{l_{\rho_1}(\beta)}{l_{\rho_3}(\beta)}, \\
\Bigg| \log \frac{l_{\rho_1}(\beta)}{l_{\rho_2}(\beta)}\Bigg|+
\Bigg| \log \frac{l_{\rho_2}(\beta)}{l_{\rho_3}(\beta)}\Bigg| &
\geq \log \frac{l_{\rho_2}(\beta)}{l_{\rho_1}(\beta)} +
\log \frac{l_{\rho_3}(\beta)}{l_{\rho_2}(\beta)}=
\log \frac{l_{\rho_3}(\beta)}{l_{\rho_1}(\beta)}.
\end{align*}

\begin{figure}[htb]
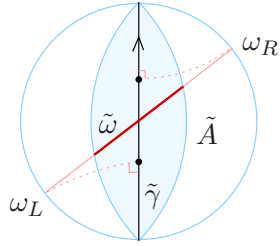

\begin{center}    
\input twistinA.pstex_t
\end{center}
\caption{Defining twist parameter.}
\label{fig:twistinA}
\end{figure}

\smallskip
Let $\gamma$ be the core of $A$ and fix a simple arc $\omega$ in $A$
that connects the two components of $\dd A$. For every $\rho \in
U_\ep(A)$, the twist parameter $tw_\rho(A)$ of $\rho$ is defined 
as follows (see also \cite[\S3]{minskyprod}). First, it is necessary
to fix an orientation of $\gamma$. Consider the universal
cover of $(A,\rho)$ in $\HH^2$ and the lifts $\tilde \gamma$,
$\tilde \omega$ of $\gamma,\omega$, respectively (see
Figure~\ref{fig:twistinA}).
Extend $\tilde \omega$ to an infinite geodesic $\bar \omega$ and let
$\omega_L, \omega_R$ be the endpoints of $\bar \omega$ that lie on the
left and right of $\tilde \gamma$, respectively.  Let $p_L,p_R$ be
respectively, the orthogonal projections of $\omega_L, \omega_R$ to
$\tilde \gamma$. Then the twist parameter is defined as
$$tw_\rho(A) = \pm \frac{d_{\HH^2}(p_L,p_R)}{l_\rho(\gamma)},$$
where the sign is $(+)$ if the direction from $p_L$ to $p_R$
coincides with the orientation of $\tilde \gamma$ and $(-)$ if
it is opposite.

Then $U_\ep(A)$ can be parameterized by the length of
$\gamma$ and the twist parameter.  The map $\rho \mapsto
(tw_\rho(A), 1/l_\rho(\gamma))$ is a homeomorphism identifying
$U_\ep(A)$ with a subset of the upper half plane:
$$U_\ep(A) = \set{ (x,y) \in \RR^2 \st y \geq \frac 1\ep }.$$

We can formulate the distance $d_{L(A)}$ on $U_\ep(A)$ in terms of these
coordinates as follows. Let $\rho_1,\rho_2 \in U_\ep(A)$ and let 
$t_i= tw_{\rho_i}(A)$, $l_i=l_{\rho_i}(\gamma)$ for $i=1,2$.
\begin{lemma}
\label{lem:annulus}
Assume that $l_1 \leq l_2$.  Then the following hold.
\renewcommand{\labelenumi}{(\roman{enumi})}
\begin{enumerate}
\item If $|t_1 - t_2 | \, l_1 \leq \log[1/l_1]$, then
$$
 d_{L(A)}(\rho_1,\rho_2) \,\add  \log \frac{l_2}{l_1}.
$$
\item If $|t_1 - t_2| \, l_1 > \log[1/l_1]$, then
$$
 d_{L(A)}(\rho_1,\rho_2) \add 
 \log \frac{ |t_1 - t_2| \, l_2}{\log[1/ l_1]}=
\log \frac{l_2}{l_1} + \log \frac{ |t_1 -t_2|\,l_1}{\log[1/l_1] }.
$$
\end{enumerate}
\end{lemma}
We remark that in comparison, the hyperbolic distance between
$z_1=(t_1,1/l_1)$ and $z_2=(t_2, 1/l_2)$ in the upper-half plane
can be estimated as follows.
Assume that $l_1 \leq l_2$.
\begin{enumerate}
\item[(i)] If $|t_1 - t_2 |\,l_1 \leq 1$ then
$$d_{\HH^2}(z_1,z_2) \add \log \frac{l_2}{l_1}.$$
\item[(ii)] If $|t_1 - t_2 |\,l_1 > 1$ then
$$d_{\HH^2}(z_1,z_2) \add \log \frac{l_2}{l_1} +
2 \log [ |t_1 - t_2| l_1].$$
\end{enumerate}
\begin{proof}
  For any arc $\beta \in {\mathcal C}(A,\,\dd A)$ intersecting
  $\gamma$ and for $\rho \in U_\ep(A)$, we can define the twist
  $tw_\rho(\beta,\gamma)$ of $\beta$ around the (oriented) curve
  $\gamma$ in the same way we defined $tw_\rho(A)$, by replacing the
  reference arc $\omega$ with $\beta$ in that construction.  It
  follows from Corollary~\ref{cor:thickthin} that
$$
\frac{l_2(\beta)}{l_1(\beta)} \mul 
\frac{\log[1/l_2] + |tw_{\rho_2}(\beta,\gamma)|\, l_2}
  {\log[1/l_1]  + |tw_{\rho_1}(\beta,\gamma)| \, l_1}.
$$
Moreover, it follows from \cite[Lemma 3.5]{minskyprod}
that 
$$| [tw_{\rho_2}(\beta,\gamma) - tw_{\rho_1}(\beta,\gamma)] -
[t_2 - t_1] | \add 1.$$
Thus, the supremum over all arcs $\beta$ intersecting $\gamma$ is
\begin{equation}
\label{eqn:supannulus}
\sup_{\beta \in
  \,{\mathcal C}(A,\,\dd A)} \frac{l_2(\beta)}{l_1(\beta)}
\mul \max \Bigg\{\frac{l_2}{l_1}, 
\frac{ \log[1/l_2]  + |t_2 - t_1| \, l_2}{\log [1/l_1]} \Bigg\}.
\end{equation}

To simplify notation, let
$$
  R_1=\frac{\log[1/l_2] + |t_2 -t_1| \, l_2}{\log[1/l_1]}, \quad
  R_2=\frac{\log[1/l_1] + |t_2 -t_1| \, l_1}{\log[1/l_2]},
$$
and
$$R=\frac{ |t_2 -t_1| \, l_2}{\log [1/l_1]}.$$
The assumption that $l_1
\leq l_2$ implies that $R_1 \mul R$ and that $R_2 < l_2/l_1 +R$.
Therefore, 
$$
  d_{L(A)} (\rho_1,\rho_2) \add \log \max \Big\{ R, \frac{l_2}{l_1}\Big\}.
$$
If $|t_1 - t_2| \, l_1 \leq \log[1/l_1]$, then $R < l_2/l_1+ 1$ so 
$$d_{L(A)} (\rho_1,\rho_2)  \add \log l_2/l_1.$$ 
If $|t_1 - t_2| \, l_1 > \log[1/l_1]$, then $R > l_2/l_1$ and hence 
\begin{equation*}
d_{L(A)} (\rho_1,\rho_2)\add \log R. \qedhere
\end{equation*}
\end{proof}

\subsection{Statement of theorem}
\label{sec:statement}
Let $\Gamma = \{\gamma_1,\ldots,\gamma_k \}$ be a collection of
disjoint, homotopically distinct simple closed curves on $S$ and let
$A_1,\dots, A_k$ be collars around $\gamma_1,\ldots,\gamma_k$,
respectively.  Choose a Fenchel-Nielsen coordinate system associated
to a marking that contains $\Gamma$ in its pants system.  Let
$s_\s(\gamma_i)$ denote the Fenchel-Nielsen twist coordinate of
$\gamma_i$. Let $U_i=U_{\ep_1}(A_i)$.  For $\s \in
Thin_{\ep_1}(S,\Gamma)$, let $\Pi_i(\s) \in U_i$ be the metric $\rho$
whose twist $tw_{\rho}(A) = s_\s(\gamma_i)$ and such that
$l_\rho(\gamma_i)=l_\s(\gamma_i)$.  Each $\s \in
Thin_{\ep_1}(S,\Gamma)$ also defines a metric
$\Pi_{S\setminus\Gamma}(\s)$ in $ \T(S\setminus\Gamma)$, obtained by
pinching the geodesic representatives of $\gamma_1,\ldots,\gamma_k$,
but otherwise leaving the metric unchanged, that is, by retaining the
same Fenchel-Nielsen coordinates.  Thus we define a homeomorphism
$$\Pi: Thin_{\ep_1}(S,\Gamma) \rightarrow \T(S\setminus\Gamma) \times
U_1 \times \cdots \times U_k.$$
Endow $\T(S
\setminus\Gamma) \times U_1 \times \cdots \times U_k$ with the sup
metric $$d_{L_\Gamma}= \sup \{ d_{L(S \setminus \Gamma)}, d_{L(A_1)},
\ldots, d_{L(A_k)} \}.$$
\begin{theorem}[product regions for Lipschitz metric]
\label{thm:prodreg}
For any $\s,\tau \in Thin_{\ep_1}(S,\Gamma)$, we have
$$d_{L}(\s,\tau) \add d_{L_\Gamma}(\Pi(\s),\Pi(\tau)).$$
\end{theorem}
The heart of the proof is Proposition~\ref{prop:main} below.

\subsection{Replacing an arc with a loop}
\begin{figure}[htb]
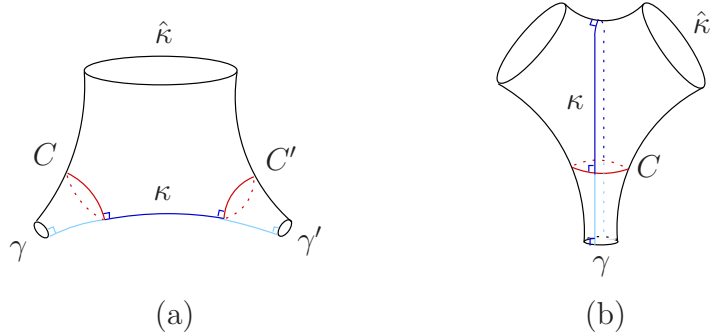

\begin{center}    
\input arctocurve.pstex_t
\end{center}
\caption{Construction of $\hat \kappa$.}
\label{fig:arctocurve}
\end{figure}

Next, we describe a procedure to replace an arc in $\zeta_Q$ with a
non-trivial, non-peripheral simple closed curve in $Q$ that has
comparable length. We assume that $Q$ is not homeomorphic to a pair of
pants.  Let $\kappa$ be a simple geodesic arc in $Q$ whose endpoints
lie in $\dd Q$ and which is perpendicular to $\dd Q$.  If the two
endpoints of $\kappa$ lie in distinct components $C,C'$ of $\dd Q$,
then the boundary of a regular neighborhood of $\kappa \cup C \cup C'$
in $Q$ consists of a single curve $\eta$. Define $\hat \kappa$ to be
the geodesic representative of $\eta$ in $S$. Note that since $Q$ is
not a pair of pants, it follows that $\eta$ is non-peripheral in $Q$,
and in particular, $\hat \kappa$ is contained in $Q$ (see
Figure~\ref{fig:arctocurve}(a)).

If both endpoints of $\kappa$ lie in a single component $C$ of $\dd
Q$, then the boundary of a regular neighborhood of $ \kappa \cup C$
has two components (see Figure~\ref{fig:arctocurve}(b)).  In this
case, define $\hat \kappa$ to be the curve of greater length between
the geodesic representatives in $S$ of the two components. Note that
$\hat \kappa$ is non-peripheral in $Q$ and in particular, it is
contained in $Q$.  Also note that unlike the preceding case, the
choice of $\hat \kappa$ depends on the geometry of the surface.

\begin{lemma}
\label{lem:arctocurve}
Suppose $Q$ is a hyperbolic component of an
$(\ep_0,\ep_1)$--decomposition of $\s$.  Let $\kappa$ be an arc in $Q$
perpendicular to $\dd Q$ and let $\hat \kappa$ be the associated
simple closed curve constructed above.  If $l_\s(\hat \kappa)>c_0$
for the Margulis constant $c_0$, 
then 
$$l_\s(\kappa) \mul l_\s(\hat \kappa).$$
\end{lemma}

\begin{proof}
Let $C,C'$ denote the components of $\dd Q$ that contain the
endpoints of $\kappa$, where we take $C=C'$ if the endpoints lie on
the same component. Let $\gamma, \gamma'$ denote the geodesic
representatives of $C,C'$ in $S$. By hypothesis, $\gamma$ and
$\gamma'$ have embedded collars in $S$, whose boundary components
each have length $\epsilon_0$.  Cut the collars in half along
$\gamma,\gamma'$ and let $\overline Q$ be the surface obtained by
attaching the half collars around $\gamma,\gamma'$ to $Q$, along
$C,C'$, respectively.  (In the case that $C \neq C'$ but $\gamma =
\gamma'$ in $S$, we attach a half-collar around $\gamma$ to each of
$C$ and $C'$.)  Since $\kappa$ intersects $\dd Q$ perpendicularly,
it has a natural extension to a (smooth) geodesic arc $\overline
\kappa$ with endpoints in $\dd \overline Q$ and perpendicular to
$\dd \overline Q$, as depicted in Figure~\ref{fig:arctocurve}.
 
First, consider the case when $C \neq C'$. Let $P$ be the pair of
pants with boundary components $\gamma, \gamma', \hat \kappa$ and
consider one of the right-angled hexagons of $P$, as in
Figure~\ref{fig:pentagon1}(a). 
\begin{figure}[htb]
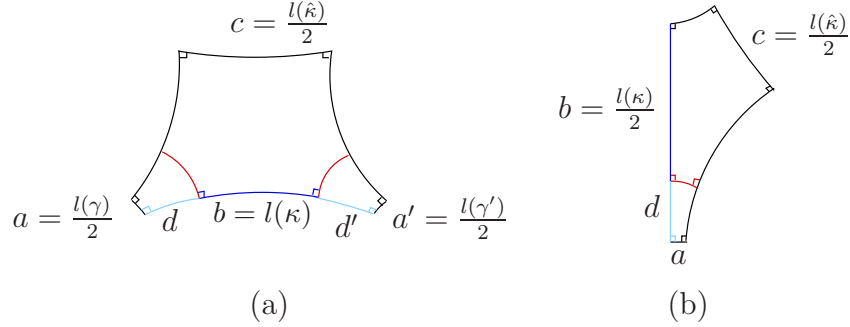

\begin{center}
\input pentagon1.pstex_t
\end{center}
\caption{A hexagon and pentagon of $P$.}
\label{fig:pentagon1}
\end{figure}

Let $a=l(\gamma)/2,a'=l(\gamma')/2$
and let $d,d'$ be the widths of the half-collars around $\gamma,
\gamma'$, respectively. Let $b=l(\kappa)$ and $c=l(\hat \kappa)/2$.
By the formula for right-angled hexagons, we have
\begin{equation}
\label{eqn:hex}
\cosh c +\cosh a \cosh a'= \sinh a \sinh a'\cosh (b+d+d').
\end{equation}
Since $a,a'<\epsilon_1/2$ and since $\ep_1$ is smaller than the
Margulis constant $c_0=0.2629 \ldots$ \cite{yamada}, we have that
$\sinh a <2a$ and $\sinh a' <2a'$. Also, by a straightforward
calculation in $\HH^2$, we have $\epsilon_0=a \cosh d = a'\cosh d'$.
Therefore, the right-hand side of Equation~\eqref{eqn:hex} satisfies
\begin{align*}
\sinh a \sinh a'\cosh (b+d+d') & > a \cdot a' \, \frac{e^{b+d+d'}}2 \\
  & > a \cdot a' \, \frac{\cosh d \cosh d'}8 \, e^b >
  \frac{\epsilon_0^2 \, e^b}8,  \\
\sinh a \sinh a'\cosh (b+d+d') & <4\, a\cdot a'\cdot e^{b+d+d'} \\
  & <a\cdot a'\cdot \cosh d \cosh d' \, e^b< \epsilon_0^2 \, e^b.
\end{align*}
On the other hand, since $a,a'< \epsilon_1/2 < \ep_0/4 <c_0/4$ and 
$c>c_0/2$, we have 
$$\cosh a \cosh a' < \cosh(a+a') < \cosh \frac{c_0}{2} < \cosh c.$$
Therefore, Equation~\eqref{eqn:hex} combined with the three equations
above gives
$$\epsilon_0^2 \, e^b/16 <\cosh c < \epsilon_0^2 \, e^b.$$
Hence,
$$|\,c-b\,|=\Big|\,\frac{l(\hat \kappa)}{2} - l(\kappa)\,\Big| <
2\log\frac{1}{\epsilon_0} + k$$
for some universal constant $k (=4\log 2)$.
Thus, if $l(\kappa)$ is sufficiently large, the additive error can be
absorbed into multiplicative constants to conclude 
$l(\hat \kappa) \mul l(\kappa)$.
If $l(\kappa)$ is not sufficiently large, then $l(\hat \kappa) \mul l(\kappa)$
holds almost tautologically, because $l(\hat \kappa)$ is bounded above by
$2l(\kappa)+2\ep_0$ and is bounded below, by assumption.
  
Next consider the case where $C=C'$.  Let $P$ be the geodesic pair of
pants in $S$ filled by $\overline \kappa \cup \gamma$. The arc
$\overline \kappa$ divides the two right-angled hexagons of $P$ into
four right-angled pentagons. It is easy to see that the two pentagons
that have edges originally contained in $\hat \kappa$ are isometric
to one another.  Let $X$ be either one of them, as in
Figure~\ref{fig:pentagon1}(b).  Let $b=l(\kappa)/2$, $c=l(\hat
\kappa)/2$ and let $d$ be the width of the half-collar around
$\gamma$. Let $a$ be the length of the edge of $X$ coming from
$\gamma$.  Now, by the formula for right-angled pentagons, we have
$$\cosh c = \sinh (b+d) \sinh a.$$
It is clear that $a \leq
l(\gamma)/2$ and by applying the pentagon formula to the pentagon
which together with $X$ makes up a hexagon of $P$, we see that our
choice of $\hat \kappa$ implies $a \geq l(\gamma)/4$.  Furthermore, as
before we have $l(\gamma)\cdot \cosh d = \epsilon_0$ and since
$l(\gamma) \leq \epsilon_1$, the assumption that
$\epsilon_0/\epsilon_1>2$ is sufficient to guarantee that $d$ is large
enough that $e^{b+d}/4 <\sinh (b+d)$ holds. And, as above,
$\epsilon_1$ is small enough that $a <\sinh a <2a$.  Therefore, we
have
\begin{align*}
\cosh c&=\sinh (b+d) \sinh a > \frac{e^b \, e^d}{4} a > 
  e^b \frac{\cosh d}{4} \cdot \frac{l(\gamma)}{4} > 
  \frac{e^b \epsilon_0}{16}.\\
\cosh c&=\sinh (b+d) \sinh a < e^b \, e^d a < 
  e^b \cdot 2\cosh d \cdot l(\gamma)/2 = e^b \epsilon_0. 
\end{align*}
Hence, 
$$
 |\,c-b\,|=\Big|\,\frac{l(\hat \kappa)}{2} - \frac{l(\kappa)}{2}\,\Big| <
 \log\frac{1}{\epsilon_0} +k
$$
for some universal constant $k(=4\log 2)$.  Thus we
conclude as before that $l(\hat \kappa) \mul l(\kappa)$.
\end{proof}

We remark that in the second case above, had we not chosen $\hat
\kappa$ to be the longer of the two components of $\dd P - \gamma$,
then the lemma would not be true. This can be easily seen by
considering the construction in reverse as follows.  Take a closed
curve $\a$ in $Q$ of moderate length and a very long arc $\beta$ with
one endpoint on $\a$ and the other on a component $C$ of $\dd Q$.
Construct a new arc $\kappa$ with both endpoints on $C$ by replacing
$\beta$ with two copies of itself very close together and by
connecting their two endpoints on $\a$ by the longer arc along $\a$.
It is easy to see that the pair of pants filled by $\kappa \cup C$ has
$\a$ as a boundary component, yet, $l(\a)/ l(\kappa)$ can be made
arbitrarily small.

\subsection{Proof of product region theorem for Lipschitz metric}
For any surface $\Sigma$, let ${\mathcal C}(\Sigma)$ be the set of
homotopy classes of non-peripheral, non-trivial simple closed curves
in $\Sigma$.  We are now ready to prove:
\begin{proposition}
\label{prop:main}
Suppose ${\mathcal P}$ is a partial $(\ep_0,\ep_1)$--decomposition for
both $\s, \tau \in \T(S)$.  Then
\begin{equation}
\sup_{\zeta \in \cs} \frac{l_\tau(\zeta)}{l_\s(\zeta)} 
\mul \max_{Q, A \in {\mathcal P}} \Bigg\{ \sup_{\a \in \,{\mathcal C}(Q)}
\frac{l_\tau(\a)}{l_\s(\a)}, \sup_{\beta \in \,{\mathcal C}(A,\,\dd A)} 
\frac{l_\tau(\beta)}{l_\s(\beta)} \Bigg\}.
\end{equation}
Moreover, when taking the
maximum, we may assume $Q$ is never a pair of pants.
\end{proposition}
\begin{proof}
  By Corollary~\ref{cor:thickthin}, for any curve $\zeta \in \cs$ and
  any $\rho \in \T(S)$ which has ${\mathcal P}$ as a partial
  $(\ep_0,\ep_1)$--decomposition, we have
$$
l_\r(\zeta) \mul \sum_{Q,A \in {\mathcal P}} [\,l_\r (\zeta_Q)+
l_\r(\zeta_A)\,],
$$
where $\zeta \cap A$ is written $\zeta_A$ for short.  Applied to
$\s,\tau$, this gives
\begin{equation}
\label{eqn:ratio0}
\frac{l_\tau(\zeta)}{l_\s(\zeta)}  \mul 
\frac{ \sum_{Q,A \in {\mathcal P}} [\,l_\tau (\zeta_Q)+ l_\tau(\zeta_A)\,]}
{\sum_{Q,A \in {\mathcal P}} [\,l_\s (\zeta_Q)+ l_\s(\zeta_A)\,]} 
\leq \max_{Q,A \in {\mathcal P}} \Bigg\{ \frac{l_\tau (\zeta_Q)}
{l_\s (\zeta_Q)},\frac{l_\tau (\zeta_A)}{l_\s (\zeta_A)} \Bigg\}.
\end{equation}

Fix $Q$ and write $\zeta_Q = \sum_i m_i \kappa_i + \sum_j n_j \lambda_j$, 
where $\kappa_i$ are arcs with endpoints on $\dd Q$ and $\lambda_j$
are non-peripheral simple closed curves contained in $Q$. Then
\begin{equation}
\label{eqn:ratio}
\frac{l_\tau (\zeta_Q)}{l_\s(\zeta_Q)} 
\mul \frac{\sum_i m_i \,l_\tau(\kappa_i) +
\sum_j n_j \,l_\tau (\lambda_j)}{\sum_i m_i \,l_\sigma(\kappa_i) +
\sum_j n_j \,l_\s (\lambda_j)} \leq
\max_{i,j} \Bigg\{ \frac{l_\tau(\kappa_i)}{l_\s(\kappa_i)}, 
\frac{l_\tau(\lambda_j)}{l_\s(\lambda_j)} \Bigg\}.
\end{equation} 
The idea is to show that for every $i$,
\begin{equation}
\label{eqn:ratioarc}
\frac{l_\tau(\kappa_i)}{l_\s(\kappa_i)}
\stackrel{.}{\prec} \sup_{\a \in {\mathcal C}(Q)} \frac{l_\tau(\a)}{l_\s(\a)},
\end{equation}
by replacing $\kappa=\kappa_i$ with the associated simple closed curve
$\hat \kappa = \hat \kappa_i$ in $Q$, as described above. In the case
that $Q$ is a pair of pants, it is not hard to see that there are
multiplicative constants depending only on $\ep_0$ such that
$l(\kappa) \mul i(\kappa,\dd Q)$ and so
$$
  \frac{l_\tau(\kappa)}{l_\s(\kappa)} \mul 
  \frac{ i(\kappa,\dd Q)}{i(\kappa,\dd Q)} \mul 1.
  $$
  Therefore, it is sufficient to prove
  Equation~\eqref{eqn:ratioarc} assuming that $Q$ is not a pair of
  pants, so that we may apply Lemma~\ref{lem:arctocurve}.
  
  Recall, that in the case the two endpoints of $\kappa$ lie in the
  same component of $\dd Q$, the choice of $\hat \kappa$ depends on
  the geometry of the surface. Let $\hat \kappa(\tau), \hat
  \kappa(\s)$ denote the curves associated to $\kappa$ for the two
  metrics $\tau,\s$, respectively. Note that by definition of $\hat
  \kappa$,
$$l_\s(\hat \kappa(\tau)) \leq l_\s(\hat \kappa(\s)).$$

Now, if $l_\tau(\hat \kappa(\tau)) > c_0$, then applying
Lemma~\ref{lem:arctocurve} and using the fact that $l(\hat \kappa)
\leq 2l(\kappa) + 2 \ep_0$ always holds, we have
$$
  \frac{l_\tau(\kappa)}{l_\s(\kappa)} \mul 
  \frac{l_\tau(\hat \kappa(\tau))}{l_\s(\kappa)} \precmul 
  \frac{l_\tau(\hat \kappa(\tau))}{l_\s(\hat \kappa(\s))} \leq 
  \frac{l_\tau(\hat \kappa(\tau))}{l_\s(\hat \kappa(\tau))} \leq 
  \sup_{\a \in {\mathcal C}(Q)} \frac{l_\tau(\a)}{l_\s(\a)}.
  $$
  If $l_\tau(\hat \kappa(\tau)) \leq c_0$, then in the
  $\tau$--metric, the three boundary curves of the geodesic pair of
  pants $P$ spanned by $\overline \kappa \cup \gamma \cup \gamma'$
  (see Lemma~\ref{lem:arctocurve} above) all have length shorter than
  $c_0$.  By using the formulas for right-angled pentagons and
  hexagons as in the proof of Lemma~\ref{lem:arctocurve}, it is easy
  to show that this implies $l_\tau(\kappa)$ is bounded above.
  Furthermore, since $\kappa$ meets $\dd Q$ and $\dd Q$ has an
  embedded regular neighborhood of some definite width depending on
  $\ep_0$, it follows that $l_\s(\kappa)$ is bounded below.  Hence,
$$
  \frac{l_\tau(\kappa)}{l_\s(\kappa)} \stackrel{.}{\prec}
  \frac{1}{l_\s(\kappa)} \stackrel{.}{\prec} 1.
$$
Since the ratio
$l_\tau(\kappa)/l_\s(\kappa)$ is bounded above,
Equation~\eqref{eqn:ratioarc} is tautologically satisfied in this case.
Thus Equation~\eqref{eqn:ratioarc} is proved.

Combined with Equations~\eqref{eqn:ratio0} and \eqref{eqn:ratio} 
we now have
$$
  \frac{l_\tau(\zeta)}{l_\s(\zeta)} \precmul 
  \max_{Q,A \in {\mathcal P}} \Bigg\{ \sup_{\a \in {\mathcal C}(Q)}
  \frac{l_\tau(\a)}{l_\s(\a)},\sup_{\beta \in \,{\mathcal C}(A,\,\dd A)} 
  \frac{l_\tau(\beta)}{l_\s(\beta)} \Bigg\},
  $$
  where by $l_\s(\beta),l_\tau(\beta)$ for $\beta \in {\mathcal
    C}(A,\,\dd A)$, we mean the length of $\beta$ in the metrics
  $\Pi_A(\s), \Pi_A(\tau) \in U_{\ep_1}(A)$, respectively, as defined
  in Section~\ref{sec:statement}.  Therefore, the supremum of the left
  hand side, taken over all $\zeta \in \cs$, is bounded by the
  quantity on the right hand side.

Finally, since ${\mathcal C}(Q) \subset \cs$, 
it is clear that for every $Q \in {\mathcal P}$
$$\sup_{\zeta \in \cs} \frac{l_\tau(\zeta)}{l_\s(\zeta)} \geq \max_{Q
  \in {\mathcal P}} \Bigg\{ \sup_{\a \in \,{\mathcal C}(Q)}
\frac{l_\tau(\a)}{l_\s(\a)} \Bigg\}.$$
To complete the proof, we will
show that there is a simple closed curve $\zeta$ such that
$$\frac{l_\tau(\zeta)}{l_\s(\zeta)} \stackrel{.}{\succ} \sup_{\beta
  \in \,{\mathcal C}(A,\,\dd A)} \frac{l_\tau(\beta)}{l_\s(\beta)}.$$
Given an annulus $A \in {\mathcal P}$, suppose that $\beta$ is a
geodesic arc in $A$ with endpoints on $\dd A$ that realizes the
supremum on the right.  Let $\a$ denote the core geodesic of $A$. If
$\a$ does not separate $S$, we can always find a non-trivial arc
$\delta$ contained in $S \setminus A$ joining the endpoints of $\beta$
whose length in $\s$ is bounded above by some constant $\ell$
depending only on $\ep_0$, and such that $\beta \cup \delta$ forms a
non-trivial simple closed curve.  Let $\zeta$ be the closed curve in
the isotopy class of $\beta \cup \delta$.  Then we have $l_\rho(\zeta)
\mul l_\rho(\beta)$ and it follows from Corollary~\ref{cor:thickthin}
that $l_\tau(\zeta) \stackrel{.}{\succ} l_\tau(\beta)$. Therefore,
$$\frac{l_\tau(\zeta)}{l_\s(\zeta)} \stackrel{.}{\succ}
\frac{l_\tau(\beta)}{l_\s(\beta)}.$$
In the case that $\a$ separates
$S$, we take an additional arc $\beta'$ in $(A,\dd A)$ disjoint from
$\beta$. Construct a simple closed curve $\zeta$ by joining the pairs
of endpoints of $\beta, \beta'$ which lie in a common component of
$\dd A$, by arcs $\delta,\delta'$ in $S \setminus A$, whose lengths in
$\s$ are uniformly bounded above.  By the same argument as before, we
can again show that
\begin{equation*}
\frac{l_\tau(\zeta)}{l_\s(\zeta)} \stackrel{.}{\succ} 
\frac{l_\tau(\beta)}{l_\s(\beta)}. \qedhere
\end{equation*}
\end{proof}

We conclude this section with the proof of Theorem~\ref{thm:prodreg}:

\begin{proof}[Proof of Theorem~\ref{thm:prodreg}]
By Proposition~\ref{prop:main}, we have that
$$d_L(\s,\tau) \add \max_{Q,A \in {\mathcal P}} \Big\{ \log \sup_{\a
  \in \,{\mathcal C}(Q)} \frac{l_\tau(\a)}{l_\s(\a)}, \ \log \sup_{\a
  \in \,{\mathcal C}(Q)} \frac{l_\s(\a)}{l_\tau(\a)}, \ 
d_{L(A)}(\s,\tau) \Big\}.$$
Therefore, to complete the proof, it would
be sufficient to show that
$$\sup_{\a \in \,{\mathcal C}(Q)} \frac{l_\tau(\a)}{l_\s(\a)} \mul
\sup_{\a \in \,{\mathcal C}(Q)} \frac{
  l_{\Pi_{S\setminus\Gamma}(\tau)}(\a)}{
  l_{\Pi_{S\setminus\Gamma}(\s)}(\a)}.$$
However, it was already shown
in \cite{minskyprod} that for $\rho \in Thin_\ep(S,\Gamma)$, the space
$(Q,\rho)$ embeds $K$--quasiconformally (in fact, biLipschitz), with
uniform $K$, in $(Q,\pi_{S\setminus\Gamma}(\rho))$.  Thus, the lengths
of curves in the two spaces are comparable and the theorem follows.
\end{proof}

\section{Comparison on a thin region} \label{sec:thin}
\noindent
We now provide an example that illustrates the discrepancy between the
Lipschitz and Teichm\"uller distances stated in the introduction:
\setcounter{introthm}{0}
\begin{introthm}
There are sequences $\s_n,\tau_n
\in \T(S)$ such that, as $n \rightarrow \infty$,
$$d_L(\s_n,\tau_n) \rightarrow 0, \ \ d_\T(\s_n,\tau_n) \rightarrow
\infty.$$
\end{introthm}
\begin{proof} Let $\s_n$ be a hyperbolic metric on $S$
  such that there is exactly one short curve $\gamma$ of length
  $l_{\s_n}(\gamma)=\ep_n$ and let $\tau_n={D_\gamma^{T_n}}(\s_n)$ be
  the metric obtained from $\s_n$ by $T_n$ Dehn twists around
  $\gamma$.  In this case, $l_{\s_n}(\gamma)=l_{\tau_n}(\gamma)=
  \ep_n$.  Set $\ep_n= e^{-P_n}$, $T_n=e^{P_n + q_n}$ and choose the
  sequences of positive integers $P_n,q_n$ so that
$$P_n \rightarrow \infty, \quad q_n \rightarrow \infty
\quad\text{and}\quad \frac{e^{q_n}}{P_n} \rightarrow 0 
\quad\text{as}\quad n\to \infty .$$  

On the one hand, it follows Theorem~\ref{thm:minskyprod} that
$$d_{\T}(\s_n,\tau_n) \add \log[T_n \ep_n]=q_n \to \infty.$$
On the
other hand, it follows from Proposition~\ref{prop:thickthin} that for
a simple closed curve $\zeta$ in $S$, we have
\begin{equation*}
\frac{l_{\tau_n}(\zeta)}{l_{\s_n}(\zeta)}=
\frac{l_{\tau_n}(\zeta_Q) + \big[ \log [\ep_0/\ep_n] +
\ep_n\cdot Tw_{\tau_n}(\zeta,\gamma)/2 + O(1) \big]
\cdot i(\zeta,\gamma)}
{l_{\s_n}(\zeta_Q) + \big[ \log [\ep_0/\ep_n] +
\ep_n \cdot Tw_{\s_n}(\zeta,\gamma)/2 + O(1)
\big]\cdot
i(\zeta,\gamma)},
\end{equation*}
where $O(1)$ represents an error that is independent of $\zeta$,
$\s_n$, $\tau_n$ and that is bounded in absolute value by some uniform
constant. Since $\s_n,\tau_n$ coincide outside of $A$, we have $
l_{\tau_n}(\zeta_Q)=l_{\s_n}(\zeta_Q)$. Therefore,
$$\sup_{\zeta}\frac{l_{\tau_n}(\zeta)}{l_{\s_n}(\zeta)} \leq \max
\Bigg\{ 1, \ \sup_{\zeta}\frac{2\log [\ep_0/\ep_n] + \ep_n\cdot
  Tw_{\tau_n}(\zeta,\gamma)+ O(1)}{2\log [\ep_0/\ep_n] + \ep_n \cdot
  Tw_{\s_n}(\zeta,\gamma)+ O(1)} \Bigg\}$$
and by the same reasoning
used to deduce Equation~(\ref{eqn:supannulus}), the supremum on the
right-hand side is equal to
$$\frac{2\log [1/\ep_n] + \ep_n\cdot T_n + O(1)}{2\log [1/\ep_n] +
  O(1)} =\frac{2P_n +e^{q_n} + O(1) }{2P_n + O(1)}. $$
Thus, we have
\begin{equation*}
\lim_{n \to \infty}
d_{L}(\s_n,\tau_n)=\lim_{n \to \infty} \log
\frac{2P_n +e^{q_n} + O(1) }{2P_n + O(1)} =0. \qedhere
\end{equation*}
\end{proof}

So far, we have seen that if $\s,\tau \in \T(S)$ are both in the thick
part then $d_L(\s,\tau) \asymp d_\T(\s,\tau)$, but that if $\s,\tau$
have a short curve in common, then the two distances are no longer
comparable. The following proposition shows that, in some sense, this
is the only way for the distances to diverge.
\begin{proposition}
\label{prop:onlyway}
If $\s,\tau \in \ts$ have no short curves in common, then
$d_L(\s,\tau) \asymp d_\T(\s,\tau)$.
\end{proposition}
\begin{proof}
  Let $\Gamma_\s$ be the set of curves whose length is less than
  $\ep_1$ at $\s$ and let $\bar \s$ be the point in the thick part of
  $\ts$ obtained from $\s$ by increasing the length of each curve in
  $\Gamma_\s$ to $\ep_1$ but otherwise leaving the metric unchanged.
  This, as usual, can be achieved by choosing a marking $\mu_\s$ of
  $S$ that contains $\Gamma_\s$ in its pants system and altering the
  associated Fenchel-Nielsen length coordinates as desired. We define
  $\bar \tau$ analogously by increasing the length of every short
  curve of $\tau$ to $\ep_1$. It follows from
  Theorem~\ref{thm:prodreg} and Lemma~\ref{lem:annulus} that
  $$
  d_L(\s,\bar \s) \asymp \log \max_{\a \in \Gamma_\s}
  \Bigg\{\frac{l_{\bar \s}(\a)}{ l_\s(\a)} \Bigg\}, \ d_L(\tau,\bar
  \tau) \asymp \log \max_{\a \in \Gamma_\tau}\Bigg\{ \frac{l_{\bar
      \tau}(\a)}{l_\tau(\a)} \Bigg\}.
  $$
  Since curves that are short in $\s$ are not short in $\tau$ and
  vice versa, the above equation implies that
\begin{equation}
\label{eqn:triangle0}
d_L(\s,\bar \s) \prec d_L(\s,\tau)
\mbox{ and } d_L(\tau,\bar \tau) \prec d_L(\s,\tau).
\end{equation}
By the triangle inequality, we also have
\begin{equation} 
\begin{split} \label{eqn:triangle}
  &d_L(\s,\tau) \geq d_L(\bar \s, \bar \tau) - d_L(\s,\bar \s) -
  d_L(\bar \tau, \tau),  \\
  &d_L(\s,\tau) \leq d_L(\bar \s, \bar \tau) + d_L(\s,\bar \s) +
  d_L(\bar \tau, \tau)
\end{split} \end{equation}
Combining Equations~\eqref{eqn:triangle0} and \eqref{eqn:triangle}, we get
\begin{equation}
\label{eqn:dL}
d_L(\s,\tau) \asymp d_L(\bar \s, \s) + d_L(\bar \s,\bar \tau)
+ d_L(\tau,\bar \tau).
\end{equation}

Analogously, it follows from Theorem~\ref{thm:minskyprod} and 
Equation~(\ref{eqn:kerckhoff}) that
$$d_\T(\s,\bar \s) \prec d_\T(\s,\tau) \mbox{ and } d_\T(\tau,\bar
\tau) \prec d_\T(\s,\tau)$$
and combined with the triangle inequality
again, we get
\begin{equation}
\label{eqn:dT}
d_\T(\s,\tau) \asymp d_\T(\bar \s, \s) + d_\T(\bar \s, \bar \tau)
+ d_\T(\tau,\bar \tau).
\end{equation}
Now, by Theorem~\ref{thm:prodreg}, Lemma~\ref{lem:annulus}, and
Theorem~\ref{thm:minskyprod} we have
$$
d_L(\s,\bar \s) \asymp d_\T(\s,\bar \s) \ \mbox{ and } \ 
d_L(\tau,\bar \tau) \asymp d_\T(\tau,\bar \tau)$$
and by Theorem~\ref{thm:thick} we have
$$d_L(\bar \s, \bar \tau) \asymp d_\T(\bar \s, \bar \tau).$$
Thus it follows from Equations~\eqref{eqn:dL} and \eqref{eqn:dT} that
$d_L(\s,\tau) \asymp d_\T(\s,\tau)$, as claimed.
\end{proof}


\end{document}

%% file: twistinA.pstex_t
\begin{picture}(0,0)%
\includegraphics{twistinA.pstex}%
\end{picture}%
\setlength{\unitlength}{3947sp}%
\begingroup\makeatletter\ifx\SetFigFont\undefined%
\gdef\SetFigFont#1#2#3#4#5{%
  \reset@font\fontsize{#1}{#2pt}%
  \fontfamily{#3}\fontseries{#4}\fontshape{#5}%
  \selectfont}%
\fi\endgroup%
\begin{picture}(1560,1524)(991,-1723)
\put(2161,-1101){\makebox(0,0)[lb]{\smash{\SetFigFont{11}{13.2}{\familydefault}{\mddefault}{\updefault}$\tilde A$}}}
\put(1551,-1031){\makebox(0,0)[lb]{\smash{\SetFigFont{11}{13.2}{\familydefault}{\mddefault}{\updefault}$\tilde \omega$}}}
\put(2453,-513){\makebox(0,0)[lb]{\smash{\SetFigFont{11}{13.2}{\familydefault}{\mddefault}{\updefault}$\omega_R$}}}
\put(991,-1526){\makebox(0,0)[lb]{\smash{\SetFigFont{11}{13.2}{\familydefault}{\mddefault}{\updefault}$\omega_L$}}}
\put(1836,-1470){\makebox(0,0)[lb]{\smash{\SetFigFont{11}{13.2}{\familydefault}{\mddefault}{\updefault}$\tilde \gamma$}}}
\end{picture}

%% file: arctocurve.pstex_t
\begin{picture}(0,0)%
\includegraphics{arctocurve.pstex}%
\end{picture}%
\setlength{\unitlength}{3947sp}%
\begingroup\makeatletter\ifx\SetFigFont\undefined%
\gdef\SetFigFont#1#2#3#4#5{%
  \reset@font\fontsize{#1}{#2pt}%
  \fontfamily{#3}\fontseries{#4}\fontshape{#5}%
  \selectfont}%
\fi\endgroup%
\begin{picture}(4375,2077)(365,-1540)
\put(1276,280){\makebox(0,0)[lb]{\smash{\SetFigFont{12}{14.4}{\familydefault}{\mddefault}{\updefault}$\hat \kappa$}}}
\put(1268,-710){\makebox(0,0)[lb]{\smash{\SetFigFont{12}{14.4}{\familydefault}{\mddefault}{\updefault}$\kappa$}}}
\put(365,-1055){\makebox(0,0)[lb]{\smash{\SetFigFont{12}{14.4}{\familydefault}{\mddefault}{\updefault}$\gamma$}}}
\put(2172,-1018){\makebox(0,0)[lb]{\smash{\SetFigFont{12}{14.4}{\familydefault}{\mddefault}{\updefault}$\gamma'$}}}
\put(3860,-140){\makebox(0,0)[lb]{\smash{\SetFigFont{12}{14.4}{\familydefault}{\mddefault}{\updefault}$\kappa$}}}
\put(4655,340){\makebox(0,0)[lb]{\smash{\SetFigFont{12}{14.4}{\familydefault}{\mddefault}{\updefault}$\hat \kappa$}}}
\put(4025,-1152){\makebox(0,0)[lb]{\smash{\SetFigFont{12}{14.4}{\familydefault}{\mddefault}{\updefault}$\gamma$}}}
\put(1276,-1486){\makebox(0,0)[lb]{\smash{\SetFigFont{12}{14.4}{\familydefault}{\mddefault}{\updefault}(a)}}}
\put(3976,-1486){\makebox(0,0)[lb]{\smash{\SetFigFont{12}{14.4}{\familydefault}{\mddefault}{\updefault}(b)}}}
\put(511,-490){\makebox(0,0)[lb]{\smash{\SetFigFont{12}{14.4}{\familydefault}{\mddefault}{\updefault}$C$}}}
\put(1970,-521){\makebox(0,0)[lb]{\smash{\SetFigFont{12}{14.4}{\familydefault}{\mddefault}{\updefault}$C'$}}}
\put(4302,-567){\makebox(0,0)[lb]{\smash{\SetFigFont{12}{14.4}{\familydefault}{\mddefault}{\updefault}$C$}}}
\end{picture}

%% file: pentagon1.pstex_t
\begin{picture}(0,0)%
\includegraphics{pentagon1.pstex}%
\end{picture}%
\setlength{\unitlength}{3947sp}%
\begingroup\makeatletter\ifx\SetFigFont\undefined%
\gdef\SetFigFont#1#2#3#4#5{%
  \reset@font\fontsize{#1}{#2pt}%
  \fontfamily{#3}\fontseries{#4}\fontshape{#5}%
  \selectfont}%
\fi\endgroup%
\begin{picture}(4793,2001)(237,-1540)
\put(1606,287){\makebox(0,0)[lb]{\smash{\SetFigFont{12}{14.4}{\familydefault}{\mddefault}{\updefault}$c=\frac{l(\hat \kappa)}{2}$}}}
\put(237,-935){\makebox(0,0)[lb]{\smash{\SetFigFont{12}{14.4}{\familydefault}{\mddefault}{\updefault}$a=\frac{l(\gamma)}{2}$}}}
\put(1175,-966){\makebox(0,0)[lb]{\smash{\SetFigFont{12}{14.4}{\familydefault}{\mddefault}{\updefault}$d$}}}
\put(2622,-943){\makebox(0,0)[lb]{\smash{\SetFigFont{12}{14.4}{\familydefault}{\mddefault}{\updefault}$a'=\frac{l(\gamma')}{2}$}}}
\put(1497,-905){\makebox(0,0)[lb]{\smash{\SetFigFont{11}{13.2}{\familydefault}{\mddefault}{\updefault}$b=l(\kappa)$}}}
\put(2262,-1002){\makebox(0,0)[lb]{\smash{\SetFigFont{12}{14.4}{\familydefault}{\mddefault}{\updefault}$d'$}}}
\put(4880,213){\makebox(0,0)[lb]{\smash{\SetFigFont{12}{14.4}{\familydefault}{\mddefault}{\updefault}$c=\frac{l(\hat \kappa)}{2}$}}}
\put(4369,-1167){\makebox(0,0)[lb]{\smash{\SetFigFont{12}{14.4}{\familydefault}{\mddefault}{\updefault}$a$}}}
\put(4204,-875){\makebox(0,0)[lb]{\smash{\SetFigFont{12}{14.4}{\familydefault}{\mddefault}{\updefault}$d$}}}
\put(3665,-253){\makebox(0,0)[lb]{\smash{\SetFigFont{12}{14.4}{\familydefault}{\mddefault}{\updefault}$b=\frac{l(\kappa)}{2}$}}}
\put(1726,-1486){\makebox(0,0)[lb]{\smash{\SetFigFont{12}{14.4}{\familydefault}{\mddefault}{\updefault}(a)}}}
\put(4351,-1486){\makebox(0,0)[lb]{\smash{\SetFigFont{12}{14.4}{\familydefault}{\mddefault}{\updefault}(b)}}}
\end{picture}